\documentclass[review]{elsarticle}

\usepackage[top=3cm,bottom=3cm,left=3cm,right=3cm]{geometry}

\usepackage{algorithm}
\usepackage[noend]{algpseudocode}
\usepackage{graphicx}
\usepackage{float}
\usepackage{amsmath,amsfonts,amsthm,amssymb}
\usepackage{tikz}
\usetikzlibrary{arrows.meta}
\usepackage{mdframed}
\journal{Discrete Mathematics}

\newtheorem{theorem}{Theorem}
\newtheorem{lemma}[theorem]{Lemma}
\newtheorem{corollary}[theorem]{Corollary}

\newtheorem{proposition}[theorem]{Proposition}

\begin{document}

\begin{frontmatter}

\title{Well-ve-Dominated Graphs}

\author[gtu]{Yasemin Büyükçolak\corref{cor1}}
\ead{y.buyukcolak@gtu.edu.tr}
\cortext[cor1]{Corresponding author}

\address[gtu]{Department of Mathematics, Gebze Technical University,
41400 Gebze, Kocaeli, T\"urkiye}

\begin{abstract}
Vertex-edge domination is a natural variant of domination in which a vertex
\emph{ve-dominates} an edge whenever it is incident to the edge or adjacent to one
of its endpoints. A set of vertices is a ve-dominating set if it ve-dominates
every edge of the graph.
In this paper, we introduce the notion of \emph{well-ve-dominated graphs},
namely graphs in which all minimal ve-dominating sets have the same
cardinality. Equivalently, these are graphs whose minimal isolating sets all
have the same cardinality; we refer to them as \emph{well-isolated graphs}. We
prove that recognizing well-ve-dominated graphs is co-NP-complete. We determine perfectly
well-ve-dominated graphs, the largest hereditary subclass, and obtain both a
complete forbidden induced subgraph characterization and a linear-time
recognition algorithm. In particular, every connected nontrivial member of this subclass
has vertex-edge domination and isolation number one. Finally, we characterize the class of 
well-ve-dominated trees. More precisely, we prove that well-ve-domination and well-ve-coveredness coincide on nontrivial trees, where well-ve-coveredness requires all minimal independent ve-dominating sets to have the same cardinality. The characterization
yields a linear-time recognition algorithm and identifies the reduced
members with the extremal trees attaining equality in the sharp isolation
bound.
\end{abstract}

\end{frontmatter}
% Discrete Mathematics requests an alphabetically ordered, numbered
% bibliography. Since elsarticle-num follows first-citation order, the
% following list establishes the required alphabetical order.
\nocite{AGRAWAL2024106419,boutrig2018,
BoutrigChellaliHaynesHedetniemi2016,BoutrigChellaliMeddah2025,
BoyerGoddard2024,Buyukcolak2025,CaroHansberg2017,
CorneilLerchsStewart1981,CorneilPerlStewart1985,
FinbowHartnellNowakowski1988,FrendrupHartnellVestergaard2010,Jena2022,
krishnakumari2017,krishnakumari2014,LemanskaEtAl2021,
LeskPlummerPulleyblank1984,lewis2007,lewis2010,Paul2023,peters1986,
Plummer1970,West2001,Ziemann2020,Zylinski2019}
\section{Introduction}\label{sec:intro}

Uniformity conditions on minimal or maximal graph-theoretic structures play
a central role in structural graph theory. Classical examples include
well-covered graphs, introduced by Plummer~\cite{Plummer1970}, in which all
maximal independent sets have the same cardinality; equimatchable
graphs~\cite{LeskPlummerPulleyblank1984}; well-dominated
graphs~\cite{FinbowHartnellNowakowski1988}; and well-edge-dominated
graphs~\cite{FrendrupHartnellVestergaard2010}. Such conditions often reveal
structural regularities and lead to interesting recognition problems.

In this paper, we study the corresponding uniformity question for
vertex-edge domination. Vertex-edge domination was introduced by
Peters~\cite{peters1986} and further developed by Lewis~\cite{lewis2007}.
A vertex $v$ \emph{ve-dominates} an edge if it is incident to the edge or
adjacent to one of its endpoints, and a set of vertices is a
\emph{ve-dominating set} if it ve-dominates every edge of the graph.
Vertex-edge domination has been investigated from extremal, structural, and
algorithmic perspectives; see, for example,
\cite{lewis2010,krishnakumari2014,BoutrigChellaliHaynesHedetniemi2016,
boutrig2018,krishnakumari2017,Ziemann2020,Jena2022,Paul2023,
Buyukcolak2025}.

The first uniformity notion for vertex-edge domination was recently
introduced by Boutrig et al.~\cite{BoutrigChellaliMeddah2025}. The minimum
cardinality of an independent ve-dominating set of a graph $G$ is denoted by
$i_{\mathrm{ve}}(G)$, while the maximum cardinality of a minimal independent
ve-dominating set is denoted by $\beta_{\mathrm{ve}}(G)$. A graph $G$ is
\emph{well-ve-covered} if all its minimal independent ve-dominating sets
have the same cardinality, or equivalently if
$i_{\mathrm{ve}}(G)=\beta_{\mathrm{ve}}(G)$. Boutrig et al. proved that
recognizing well-ve-covered graphs is $\mathrm{co\text{-}NP}$-complete and
obtained a constructive characterization of well-ve-covered trees. However,
minimal ve-dominating sets need not be independent. Thus, well-ve-coveredness
does not address the stronger uniformity question of when \emph{all} minimal
ve-dominating sets have the same cardinality.

Motivated by this question, we introduce the notion of \emph{well-ve-domination}. We call a graph $G$ \emph{well-ve-dominated} if
all its minimal ve-dominating sets have the same cardinality, or equivalently
if $\gamma_{\mathrm{ve}}(G)=\Gamma_{\mathrm{ve}}(G)$, where
$\gamma_{\mathrm{ve}}(G)$ and $\Gamma_{\mathrm{ve}}(G)$ denote,
respectively, the minimum and maximum cardinalities of minimal ve-dominating
sets. Since $\gamma_{\mathrm{ve}}(G)\leq i_{\mathrm{ve}}(G)
\leq\beta_{\mathrm{ve}}(G)\leq\Gamma_{\mathrm{ve}}(G)$, every
well-ve-dominated graph is indeed well-ve-covered. However, the converse is false. For
example, the graph $G=\overline{P_3\cup P_3}$ in
Figure~\ref{fig:well-ve-covered-not-dominated} is well-ve-covered because
all its minimal independent ve-dominating sets have cardinality one, whereas
$\{a\}$ and $\{c,f\}$ are minimal ve-dominating sets of different
cardinalities. Hence, well-ve-dominated graphs form a proper subclass of
well-ve-covered graphs.

\begin{figure}[b]
\centering
\resizebox{0.45\linewidth}{!}{%
\begin{tikzpicture}[
    vertexL/.style={
        circle,
        draw,
        fill=green!12,
        minimum size=5mm,
        inner sep=0pt
    },
    vertexR/.style={
        circle,
        draw,
        fill=violet!12,
        minimum size=5mm,
        inner sep=0pt
    },
    edge/.style={
        black,
        line width=0.9pt
    },
    nonedge/.style={
        dashed,
        gray!70,
        line width=1pt
    },
    every node/.style={font=\small}
]

\node[vertexL] (a) at (0,1.5) {$a$};
\node[vertexL] (b) at (0,-1.5) {$b$};
\node[vertexL] (c) at (-2.2,0) {$c$};

\node[vertexR] (d) at (4,1.5) {$d$};
\node[vertexR] (e) at (4,-1.5) {$e$};
\node[vertexR] (f) at (6.2,0) {$f$};

\foreach \x/\y in {a/b,a/d,a/e,a/f,
                   b/d,b/e,b/f,
                   c/d,c/e,c/f,
                   d/e}
{
    \draw[edge] (\x)--(\y);
}

\end{tikzpicture}
}
\caption{The graph $\overline{P_3\cup P_3}$ is well-ve-covered but not
well-ve-dominated. }
\label{fig:well-ve-covered-not-dominated}
\end{figure}

On the other hand, vertex-edge domination is a special case of graph isolation in which
the forbidden family consists of a single edge. More precisely, a set
$D\subseteq V(G)$ is an isolating set if the vertices outside its closed
neighborhood induce no edges, which is equivalent to $D$ being a
ve-dominating set. We therefore use the term \emph{well-isolated} for a
graph whose minimal isolating sets all have the same cardinality. Thus,
well-isolation is the isolation-theoretic formulation of well-ve-domination,
rather than a distinct graph class.

The literature on graph isolation has focused primarily on the minimum
cardinality of isolating sets and on extremal bounds for the isolation
number. To the best of our knowledge, uniformity among all minimal isolating
sets has not previously been investigated. Thus,
well-ve-domination also connects the uniformity viewpoint familiar from
well-covered and well-dominated graphs with the extremal and structural
theory of graph isolation.

From the extremal point of view,
Caro and
Hansberg~\cite{CaroHansberg2017} proved that every connected graph $G$ of
order $n\geq3$, other than $C_5$, satisfies $\iota(G)\leq n/3$, where
$\iota(G)$ denotes the minimum cardinality of an isolating set of $G$.
Żyliński~\cite{Zylinski2019} independently established the same bound for
connected graphs of order at least six in the language of vertex-edge
domination. Furthermore, the trees attaining equality were characterized in
\cite{krishnakumari2014,LemanskaEtAl2021}: they are obtained from an
arbitrary tree by attaching a pendant path of length two to each of its
vertices. This extremal family will reappear naturally in our
characterization of reduced well-ve-dominated trees.

Our contributions are as follows. We introduce well-ve-dominated graphs,
establish basic structural properties, and show that the class is
not hereditary. We also prove that recognizing well-ve-dominated graphs is
$\mathrm{co\text{-}NP}$-complete. Moreover,
we determine the largest hereditary subclass of well-ve-dominated
graphs, namely \emph{perfectly well-ve-dominated graphs}. We characterize
this subclass by three forbidden induced subgraphs and give an equivalent
cotree characterization,
which yields a linear-time recognition algorithm. In particular, although
recognizing well-ve-dominated graphs is $\mathrm{co\text{-}NP}$-complete in
general, its largest hereditary subclass admits a rigid structural
description and efficient recognition. Moreover, every connected nontrivial
member of this subclass satisfies $\gamma_{\mathrm{ve}}=\iota=1$, and its
triangle-free members are precisely the complete bipartite graphs.
Finally, building on the characterization of well-ve-covered trees by
Boutrig et al.~\cite{BoutrigChellaliMeddah2025}, we prove that
well-ve-coveredness and well-ve-domination coincide on nontrivial trees. Particularly, we present a canonical representation of reduced
well-ve-dominated trees and a linear-time recognition algorithm. The
canonical form identifies the reduced well-ve-dominated trees of order at
least three precisely with the extremal family satisfying
$\iota(T)=|V(T)|/3$. Thus, in the reduced tree setting, uniformity among
minimal isolating sets coincides with equality in the sharp isolation bound.

The rest of the paper is organized as follows.
Section~\ref{sec:pre} introduces the notation, basic definitions, and
structural properties used throughout the paper.
Section~\ref{sec:complexity} establishes the complexity of recognizing
well-ve-dominated graphs. Section~\ref{sec:perfectly} characterizes
perfectly well-ve-dominated graphs and gives the corresponding recognition
algorithm. Finally, Section~\ref{sec:WEDtrees} characterizes
well-ve-dominated trees and develops their canonical and algorithmic
structure.

\section{Preliminaries}\label{sec:pre}

All graphs considered in this paper are finite, simple, and undirected. For
a graph $G$, its vertex set and edge set are denoted by $V(G)$ and $E(G)$,
respectively, and its order is $|V(G)|$. The open and closed neighborhoods
of a vertex $v$ are denoted by $N_G(v)$ and $N_G[v]$, respectively. For
$X\subseteq V(G)$, let
$N_G(X)=\bigcup_{x\in X}N_G(x)$ and $N_G[X]=N_G(X)\cup X$.
The degree of $v$ is denoted by $d_G(v)$. When the underlying graph is
clear, we omit the subscript $G$. For $X\subseteq V(G)$, the subgraph
induced by $X$ is denoted by $G[X]$, and $G-X$ is obtained by deleting the
vertices of $X$ and all incident edges. We write $G-v$ for $G-\{v\}$.
The complement of $G$ is denoted by $\overline{G}$.
For vertex-disjoint
graphs $G$ and $H$, their disjoint union is denoted by $G\cup H$. We use
$P_n$, $C_n$, $K_n$, and $K_{m,n}$ for the path, cycle, complete graph,
and complete bipartite graph of the indicated orders. We write $2K_2$ for
the disjoint union of two copies of $K_2$. For vertex-disjoint graphs $G$
and $H$, their \emph{join}, denoted by $G\vee H$, is obtained from $G\cup H$
by adding every edge between $V(G)$ and $V(H)$. 
A vertex of degree one is a \emph{leaf}, and a vertex adjacent to a leaf is
a \emph{support vertex}. An edge incident to a leaf is a \emph{pendant
edge}, and it is called \emph{good} if its non-leaf endpoint has degree two.

A graph is \emph{$H$-free} if it has no induced subgraph isomorphic to $H$,
and it is \emph{$(H_1,\ldots,H_k)$-free} if it is $H_i$-free for every
$i\in\{1,\ldots,k\}$. A \emph{cograph} is a $P_4$-free graph, and every
connected cograph of order at least two is the join of two nonempty induced
subgraphs~\cite{CorneilLerchsStewart1981}. A graph is \emph{chordal} if it
has no induced cycle of length at least four. A vertex is \emph{simplicial}
if its neighborhood is a clique; every nonempty chordal graph has a simplicial
vertex~\cite{West2001}.
A \emph{cotree} of a cograph $G$ is a rooted tree whose leaves are the
vertices of $G$ and whose internal nodes are labeled by disjoint union or
join; internal nodes with these labels are called \emph{union nodes} and
\emph{join nodes}, respectively. For a cotree node $u$, let $G_u$ denote the
subgraph of $G$ induced by the leaves in the subtree rooted at $u$.
According to the label of $u$, the
graph $G_u$ is the disjoint union or the join of the graphs associated with
the children of $u$. A \emph{reduced cotree} is one in which every internal
node has at least two children and adjacent internal nodes have different
labels. A cograph can be recognized and a reduced cotree can be constructed
in $O(n+m)$ time~\cite{CorneilPerlStewart1985}.

A vertex $v$ \emph{ve-dominates} an edge $e$ if at least one endpoint of
$e$ belongs to $N_G[v]$. A set $D\subseteq V(G)$ is a
\emph{ve-dominating set} if every edge of $G$ is ve-dominated by a vertex
of $D$. The minimum cardinality of such a set is the
\emph{vertex-edge domination number} $\gamma_{\mathrm{ve}}(G)$.
A ve-dominating set is \emph{minimal} if none of its proper subsets is
ve-dominating, and the maximum cardinality of a minimal ve-dominating set
is denoted by $\Gamma_{\mathrm{ve}}(G)$.
For a ve-dominating set $S$, an edge $e$ is a \emph{private edge} of
$v\in S$ with respect to $S$ if $v$ ve-dominates $e$ and no vertex of
$S\setminus\{v\}$ ve-dominates $e$. 

\begin{proposition}[{\cite{BoutrigChellaliHaynesHedetniemi2016}}]
\label{prop:minimal-ved}
A ve-dominating set $S$ of a graph $G$ is minimal if and only if every
vertex of $S$ has a private edge with respect to $S$.
\end{proposition}

\noindent
A ve-dominating set is \emph{independent} if its vertices are pairwise
nonadjacent. The minimum cardinality of an independent ve-dominating set is
denoted by $i_{\mathrm{ve}}(G)$, and the maximum cardinality of a minimal
independent ve-dominating set is denoted by $\beta_{\mathrm{ve}}(G)$. Following \cite{BoutrigChellaliMeddah2025}, a graph is called \emph{well-ve-covered} if all its minimal independent
ve-dominating sets have the same cardinality. Equivalently,
$i_{\mathrm{ve}}(G)=\beta_{\mathrm{ve}}(G)$. For every graph $G$,
$\gamma_{\mathrm{ve}}(G)\le i_{\mathrm{ve}}(G)\le
\beta_{\mathrm{ve}}(G)\le \Gamma_{\mathrm{ve}}(G)$.
In this work, we introduce the notion of well-ve-domination.
We call a graph $G$ \emph{well-ve-dominated} if all its minimal
ve-dominating sets have the same cardinality. Thus, $G$ is
well-ve-dominated if and only if
$\gamma_{\mathrm{ve}}(G)=\Gamma_{\mathrm{ve}}(G)$.
Figure \ref{fig:well-ve-covered-not-dominated} shows that the class of well-ve-dominated graphs is a \emph{proper} subclass of well-ve-covered graphs.
\noindent

Let $\mathcal{F}$ be a family of graphs. A set $D\subseteq V(G)$ is an
\emph{$\mathcal{F}$-isolating set} if $G-N_G[D]$ contains no member of
$\mathcal{F}$ as a subgraph. The minimum cardinality of such a set is the
\emph{$\mathcal{F}$-isolation number} $\iota(G,\mathcal{F})$. When
$\mathcal{F}=\{K_2\}$, an $\mathcal{F}$-isolating set is simply called an
\emph{isolating set}, and the minimum cardinality of such a set is denoted by
$\iota(G)$.
A set $D$ is isolating if and only if $G-N_G[D]$ is edgeless. This is
equivalent to $D$ ve-dominating every edge of $G$. Hence,
$\gamma_{\mathrm{ve}}(G)=\iota(G)$.
Moreover, the minimal isolating sets of $G$ are precisely its minimal
ve-dominating sets. We call a graph \emph{well-isolated} if all its minimal
isolating sets have the same cardinality. Consequently, a graph is
well-isolated if and only if it is well-ve-dominated.

\medskip

We now proceed by establishing some basic structural properties of well-ve-dominated
graphs that will be used throughout the paper. 

\begin{proposition}\label{prop:components}
A graph $G$ is well-ve-dominated if and only if each nontrivial connected
component of $G$ is well-ve-dominated.
\end{proposition}
\begin{proof}
Let $G_1,\ldots,G_k$ be the nontrivial connected components of $G$.
The minimal ve-dominating sets of $G$ are precisely the sets
\(
D=D_1\cup\cdots\cup D_k,
\)
where each $D_i$ is a minimal ve-dominating set of $G_i$. Therefore, all
minimal ve-dominating sets of $G$ have the same cardinality if and only if
this holds for each $G_i$. Indeed, if some component $G_j$ has two minimal
ve-dominating sets of different cardinalities, then fixing a minimal
ve-dominating set in every other nontrivial component produces two minimal
ve-dominating sets of different cardinalities in $G$.
\end{proof}

Two distinct vertices $x$ and $y$ are \emph{false twins} if
$N_G(x)=N_G(y)$. A graph is \emph{reduced} if it has no false twins. Its
\emph{reduction} $R(G)$ is obtained by retaining exactly one representative
from each equivalence class of vertices with identical open neighborhoods.
It is uniquely determined up to isomorphism.

\begin{lemma}\label{lem:reduced}
A graph $G$ is well-ve-dominated if and only if its reduction $R(G)$ is
well-ve-dominated.
\end{lemma}

\begin{proof}
It suffices to show that deleting one of two false twins preserves the
possible cardinalities of minimal ve-dominating sets. Let $x$ and $y$ be
false twins and put $G'=G-y$. Since $N_G(x)=N_G(y)$, the vertices $x$ and
$y$ are nonadjacent and ve-dominate the same edges.
For every $D\subseteq V(G')$, the set $D$ ve-dominates $G'$ if and only if
it ve-dominates $G$. Indeed, every edge of $G$ not in $G'$ has the form
$yz$, with $xz\in E(G')$, and any vertex of $D$ that ve-dominates $xz$
also ve-dominates $yz$.
Now let $S$ be a minimal ve-dominating set of $G$. Since $x$ and $y$
ve-dominate the same edges, $S$ cannot contain both. If $y\in S$, replace
$y$ by $x$; otherwise leave $S$ unchanged. The resulting set
$S'\subseteq V(G')$ has $|S'|=|S|$. Since $x$ and $y$ ve-dominate exactly
the same edges, $S'$ ve-dominates $G$. It is also minimal in $G$. Indeed,
suppose that a proper subset $A\subsetneq S'$ ve-dominates $G$. If
$x\notin A$, then $A\subsetneq S$, contradicting the minimality of $S$. If
$x\in A$ and $y\in S$, then replacing $x$ by $y$ in $A$ produces a proper
subset of $S$ that ve-dominates the same edges, again a
contradiction. In the remaining case $S'=S$, and the contradiction is
immediate.
Since $S'\subseteq V(G')$, the equivalence above now implies that $S'$ is a
minimal ve-dominating set of $G'$. Conversely, if a minimal ve-dominating
set of $G'$ had a proper subset that ve-dominated $G$, then the equivalence
above would imply that the same subset ve-dominates $G'$. Hence, every
minimal ve-dominating set of $G'$ is also minimal in $G$.
Therefore, $G$ and $G'$ have the same possible cardinalities of minimal
ve-dominating sets. Repeating this deletion within each false-twin class
proves the result.
\end{proof}

\begin{proposition}\label{prop:not-hereditary}
The class of well-ve-dominated graphs, or equivalently, the class of
well-isolated graphs, is not hereditary.
\end{proposition}

\begin{proof}
The path $P_6$ is well-ve-dominated, since all its minimal ve-dominating
sets have cardinality two. 
In contrast, if $P_4=v_1v_2v_3v_4$, then
$\{v_2\}$ and $\{v_1,v_4\}$ are minimal ve-dominating sets of cardinalities
one and two, respectively. Since $P_4$
is an induced subgraph of $P_6$, the result follows.
\end{proof}

\section{Complexity of Well-ve-Dominated Graphs}\label{sec:complexity}

In this section, we investigate the computational complexity of recognizing
well-ve-dominated graphs. Equivalently, since well-ve-domination and
well-isolation coincide, we study the complexity of recognizing
well-isolated graphs.

We first recall the additional graph-theoretic terminology used in this
section. A set $D\subseteq V(G)$ is a \emph{dominating set} of $G$ if
$N_G[D]=V(G)$. A dominating set is \emph{minimal} if none of its proper
subsets is a dominating set. A graph is \emph{well-dominated} if all its
minimal dominating sets have the same cardinality.
The \emph{lexicographic product} of two graphs $G$ and $F$, denoted by
$G\circ F$, is the graph with vertex set $V(G)\times V(F)$ in which two
vertices $(u,x)$ and $(v,y)$ are adjacent if and only if either
$uv\in E(G)$, or $u=v$ and $xy\in E(F)$.

\medskip

We consider the following decision problem.

\begin{mdframed}[nobreak=true]
\noindent\textsc{Well-VE-Dominated Recognition}

\medskip

\textit{Instance:} A graph $G$.

\smallskip

\textit{Question:} Do all minimal ve-dominating sets of $G$ have the same
cardinality? Equivalently, do all minimal isolating sets of $G$ have the
same cardinality?
\end{mdframed}

The recognition problem for well-ve-dominated graphs is closely connected to
that of well-ve-covered graphs, which was proved to be
$\mathrm{co\text{-}NP}$-complete by Boutrig et
al.~\cite{BoutrigChellaliMeddah2025}. Since every well-ve-dominated graph is
necessarily well-ve-covered, it is natural to ask whether the stronger
uniformity requirement imposed by well-ve-domination makes the corresponding
recognition problem computationally easier. The following theorem shows that
this is not the case.
\begin{theorem}\label{thm:well-ve-dominated-complexity}
Recognizing whether a graph is well-ve-dominated is
$\mathrm{co\text{-}NP}$-complete. Equivalently, recognizing whether a graph
is well-isolated is $\mathrm{co\text{-}NP}$-complete.
\end{theorem}

\begin{proof}
We first show that \textsc{Well-VE-Dominated Recognition} belongs to
$\mathrm{co\text{-}NP}$. A certificate that a graph is not
well-ve-dominated consists of two minimal ve-dominating sets of different
cardinalities. Given such a certificate, one can verify in polynomial time
that both sets are ve-dominating, that each is minimal by checking that the
deletion of any of its vertices leaves some edge non-ve-dominated, and that
their cardinalities differ. Hence, the complement of the problem belongs to
$\mathrm{NP}$, and therefore \textsc{Well-VE-Dominated Recognition} belongs
to $\mathrm{co\text{-}NP}$.

To establish $\mathrm{co\text{-}NP}$-hardness, we give a polynomial-time
reduction from \textsc{Well-Dominated Recognition}, the problem of deciding
whether all minimal dominating sets of a graph have the same cardinality,
which is known to be
$\mathrm{co\text{-}NP}$-complete~\cite{AGRAWAL2024106419}.

Let $G$ be an arbitrary graph and let $H=G\circ K_2$. For each
$v\in V(G)$, let $X_v=\{v^0,v^1\}$ denote the copy of $K_2$ corresponding
to $v$. For every edge $uv\in E(G)$, all four edges between $X_u$ and
$X_v$ belong to $H$. Thus, $H$ is obtained by replacing each vertex of $G$
with a copy of $K_2$. This construction is illustrated in
Figure~\ref{fig:C4-lexicographic-product} for $G=C_4$. For any $S\subseteq V(H)$, define its projection onto $G$ by
$\pi(S)=\{v\in V(G):S\cap X_v\neq\emptyset\}$.

\begin{figure}[t]
\centering
\resizebox{0.8\linewidth}{!}{%
\begin{tikzpicture}[
    vertexG/.style={circle, draw, fill=green!12, minimum size=5mm, inner sep=0pt},
    vertexH/.style={circle, draw, fill=violet!12, minimum size=5mm, inner sep=0pt},
    edge/.style={black, line width=0.9pt},
    transformation/.style={-{Stealth[length=3mm,width=2.2mm]}, black, line width=1pt},
    every node/.style={font=\small}
]

% Original graph G=C_4
\begin{scope}[shift={(0,0)}, scale=0.8, transform shape]
\node[vertexG] (v1) at (-1.2,1.2) {$v_1$};
\node[vertexG] (v2) at ( 1.2,1.2) {$v_2$};
\node[vertexG] (v3) at ( 1.2,-1.2) {$v_3$};
\node[vertexG] (v4) at (-1.2,-1.2) {$v_4$};

\draw[edge] (v1)--(v2);
\draw[edge] (v2)--(v3);
\draw[edge] (v3)--(v4);
\draw[edge] (v4)--(v1);

\node at (0,-2) {$G=C_4$};
\end{scope}

% Transformation arrow
\draw[transformation] (1.9,0)--(3.9,0);

% Product H=C_4\circ K_2
\begin{scope}[shift={(7.1,0)}, scale=0.8, transform shape]
\node[vertexH] (v10) at (-1.25, 2.25) {$v_1^0$};
\node[vertexH] (v11) at (-2.25, 1.25) {$v_1^1$};
\node[vertexH] (v20) at ( 1.25, 2.25) {$v_2^0$};
\node[vertexH] (v21) at ( 2.25, 1.25) {$v_2^1$};
\node[vertexH] (v30) at ( 2.25,-1.25) {$v_3^0$};
\node[vertexH] (v31) at ( 1.25,-2.25) {$v_3^1$};
\node[vertexH] (v40) at (-1.25,-2.25) {$v_4^0$};
\node[vertexH] (v41) at (-2.25,-1.25) {$v_4^1$};

\draw[edge] (v10)--(v11);
\draw[edge] (v20)--(v21);
\draw[edge] (v30)--(v31);
\draw[edge] (v40)--(v41);

\draw[edge] (v10)--(v20);
\draw[edge] (v10)--(v21);
\draw[edge] (v11)--(v20);
\draw[edge] (v11)--(v21);

\draw[edge] (v20)--(v30);
\draw[edge] (v20)--(v31);
\draw[edge] (v21)--(v30);
\draw[edge] (v21)--(v31);

\draw[edge] (v30)--(v40);
\draw[edge] (v30)--(v41);
\draw[edge] (v31)--(v40);
\draw[edge] (v31)--(v41);

\draw[edge] (v40)--(v10);
\draw[edge] (v40)--(v11);
\draw[edge] (v41)--(v10);
\draw[edge] (v41)--(v11);

\node at (-2.45, 2.35) {$X_{v_1}$};
\node at ( 2.45, 2.35) {$X_{v_2}$};
\node at ( 2.45,-2.35) {$X_{v_3}$};
\node at (-2.45,-2.35) {$X_{v_4}$};

\node at (0,-3.15) {$H=C_4\circ K_2$};
\end{scope}

\end{tikzpicture}
}
\caption{The lexicographic product $H=C_4\circ K_2$.}
\label{fig:C4-lexicographic-product}
\end{figure}

\medskip
\noindent\textit{Claim 1.}
\textit{A set $S\subseteq V(H)$ is a ve-dominating set of $H$ if and only if
$\pi(S)$ is a dominating set of $G$.}

\smallskip
\noindent\textit{Proof of Claim 1.}
Suppose first that $S$ is a ve-dominating set of $H$. For each
$v\in V(G)$, consider the edge $v^0v^1$. Since this edge is ve-dominated
by $S$, some vertex of $S$ lies in $X_v$ or in $X_u$ for some
$u\in N_G(v)$. Hence, $\pi(S)\cap N_G[v]\neq\emptyset$ for every
$v\in V(G)$, and therefore $\pi(S)$ dominates $G$.
Conversely, suppose that $\pi(S)$ dominates $G$, and let $e\in E(H)$.
Choose an endpoint of $e$ belonging to some $X_v$. Since $\pi(S)$ dominates
$G$, there exists $u\in\pi(S)\cap N_G[v]$; choose $s\in S\cap X_u$.
If $u=v$, then $s$ is incident to or adjacent to the chosen endpoint of
$e$. If $u\in N_G(v)$, then every vertex of $X_u$ is adjacent to every
vertex of $X_v$, so $s$ is adjacent to that endpoint. In either case, $s$
ve-dominates $e$. Hence, $S$ is a ve-dominating set of $H$.
\hfill$\square$

\medskip
\noindent\textit{Claim 2.}
\textit{A set $S\subseteq V(H)$ is a minimal ve-dominating set of $H$ if
and only if $\pi(S)$ is a minimal dominating set of $G$ and
$|S\cap X_v|=1$ for every $v\in\pi(S)$.}

\smallskip
\noindent\textit{Proof of Claim 2.}
Suppose first that $S$ is a minimal ve-dominating set of $H$. If
$|S\cap X_v|\geq2$ for some $v\in\pi(S)$, then deleting one vertex of
$S\cap X_v$ leaves the projection unchanged. By Claim~1, the resulting
proper subset of $S$ still ve-dominates $H$, a contradiction. Thus,
$|S\cap X_v|=1$ for every $v\in\pi(S)$.
We next show that $\pi(S)$ is minimal. Otherwise, for some
$v\in\pi(S)$, the set $\pi(S)\setminus\{v\}$ still dominates $G$. Let $s$
be the unique vertex of $S\cap X_v$. Then
$\pi(S\setminus\{s\})=\pi(S)\setminus\{v\}$, and Claim~1 implies that
$S\setminus\{s\}$ still ve-dominates $H$, contradicting the minimality of
$S$. Thus, $\pi(S)$ is a minimal dominating set of $G$.

Conversely, suppose that $\pi(S)$ is a minimal dominating set of $G$ and
$|S\cap X_v|=1$ for every $v\in\pi(S)$. By Claim~1, $S$ ve-dominates
$H$. For any $s\in S$, let $v$ be the unique vertex with $s\in X_v$.
Then $\pi(S\setminus\{s\})=\pi(S)\setminus\{v\}$, which does not dominate
$G$ by the minimality of $\pi(S)$. Claim~1 therefore implies that
$S\setminus\{s\}$ does not ve-dominate $H$. Hence, $S$ is minimal.
\hfill$\square$

\medskip

By Claim~2, every minimal ve-dominating set $S$ of $H$ contains exactly one
vertex from each $X_v$ with $v\in\pi(S)$, and therefore
$|S|=|\pi(S)|$. Conversely, every minimal dominating set $D$ of $G$ can be
lifted to a minimal ve-dominating set of $H$ by choosing one vertex from
each $X_v$ with $v\in D$. Thus, the possible cardinalities of minimal
dominating sets of $G$ are exactly the possible cardinalities of minimal
ve-dominating sets of $H$. Consequently, $G$ is well-dominated if and only
if $H=G\circ K_2$ is well-ve-dominated.
Furthermore, 
the graph $G\circ K_2$ can clearly be constructed in polynomial time.
Hence, this is a polynomial-time reduction from
\textsc{Well-Dominated Recognition} to
\textsc{Well-VE-Dominated Recognition}, proving
$\mathrm{co\text{-}NP}$-hardness. Together with membership in
$\mathrm{co\text{-}NP}$, this proves that
\textsc{Well-VE-Dominated Recognition} is
$\mathrm{co\text{-}NP}$-complete.
\end{proof}

Theorem~\ref{thm:well-ve-dominated-complexity} shows that recognizing
well-ve-dominated graphs, or equivalently well-isolated graphs, is
computationally intractable in general unless $\mathrm{P}=\mathrm{NP}$.
It places this problem alongside the recognition problems for well-dominated
graphs~\cite{AGRAWAL2024106419} and well-ve-covered
graphs~\cite{BoutrigChellaliMeddah2025}, both of which are
$\mathrm{co\text{-}NP}$-complete. Thus, extending uniformity from minimal
independent ve-dominating sets to all minimal ve-dominating sets does not
reduce the complexity of the recognition problem.

\section{Perfectly Well-ve-Dominated Graphs}\label{sec:perfectly}

By Proposition~\ref{prop:not-hereditary}, the class of well-ve-dominated
graphs is not hereditary.
This naturally raises the question of determining its largest hereditary
subclass. This question is particularly relevant in view of
Theorem~\ref{thm:well-ve-dominated-complexity}:
although recognizing well-ve-dominated graphs is co-NP-complete in general,
the hereditary restriction turns out to admit a remarkably rigid structure,
including a finite forbidden induced subgraph characterization and a
linear-time recognition algorithm.

A graph $G$ is called \emph{perfectly well-ve-dominated} if every induced
subgraph of $G$ is well-ve-dominated. Analogously, $G$ is called
\emph{perfectly well-isolated} if every induced subgraph of $G$ is
well-isolated. Since well-ve-domination and well-isolation are equivalent,
these two notions coincide. Accordingly, we state the subsequent results
in terms of perfectly well-ve-dominated graphs; their perfectly well-isolated
formulations follow immediately. By definition, perfectly well-ve-dominated
graphs form the largest hereditary subclass of well-ve-dominated graphs.

As observed in the proof of Proposition~\ref{prop:not-hereditary}, $P_4$ is
not well-ve-dominated which immediately yields the
following necessary condition.

\begin{corollary}\label{cor:perfect-p4free}
Every perfectly well-ve-dominated graph is $P_4$-free. Hence, every
such graph is a cograph.
\end{corollary}

\noindent
The condition in Corollary~\ref{cor:perfect-p4free} is not sufficient.
Consider the butterfly (or bow tie) graph
$K_1\vee 2K_2$. If $x$ is its universal vertex and $ab$ and $cd$ are the
two edges of $2K_2$, then $\{x\}$ and $\{a,c\}$ are minimal ve-dominating
sets of different cardinalities. Thus, $K_1\vee 2K_2$ is $P_4$-free but not
well-ve-dominated. To see that it is $P_4$-free, observe that an induced
$P_4$ cannot contain the universal vertex, while the remaining graph
$2K_2$ contains no induced $P_4$.
Consequently, every perfectly well-ve-dominated graph is
$(P_4,K_1\vee 2K_2)$-free.
This stronger condition is still not sufficient: the graph
$\overline{P_3\cup P_3}$ is
$(P_4,K_1\vee 2K_2)$-free. Since $P_4$ is self-complementary, an induced
$P_4$ in $\overline{P_3\cup P_3}$ would correspond to an induced $P_4$ in
$P_3\cup P_3$, which is impossible. Moreover,
$\overline{K_1\vee 2K_2}=K_1\cup C_4$, and $P_3\cup P_3$ contains no induced
copy of $K_1\cup C_4$. As shown in
Figure~\ref{fig:well-ve-covered-not-dominated}, it has minimal
ve-dominating sets $\{a\}$ and $\{c,f\}$ of different cardinalities and is
therefore not well-ve-dominated. Hence, we obtain the following stronger
necessary condition.

\begin{corollary}\label{cor:perfectly-forbidden}
Every perfectly well-ve-dominated graph is
$(P_4,K_1\vee 2K_2,\overline{P_3\cup P_3})$-free.
\end{corollary}

\noindent
We will next prove that these three graphs form the complete set of forbidden
induced subgraphs. We begin with an auxiliary lemma.

\begin{lemma}\label{lem:p4-2k2-free-vedom-vertex}
Every connected $(P_4,2K_2)$-free graph has a ve-dominating vertex.
\end{lemma}

\begin{proof}
Let $G$ be a connected $(P_4,2K_2)$-free graph and let $H=\overline{G}$.
Since $P_4$ is self-complementary and
$\overline{2K_2}=C_4$, the graph $H$ is $P_4$-free and $C_4$-free.
Every induced cycle of length at least five contains an induced $P_4$.
Hence, $H$ has no induced cycle of length at least four and is chordal. It
therefore has a simplicial vertex $v$, so $N_H(v)$ is a clique in $H$. Since
\(
N_H(v)=V(G)\setminus N_G[v],
\)
the set $V(G)\setminus N_G[v]$ is independent in $G$. Hence every edge of
$G$ has at least one endpoint in $N_G[v]$. By definition, $v$
ve-dominates every edge of $G$, and thus $v$ is a ve-dominating vertex.
\end{proof}

We can now prove the forbidden induced subgraph characterization. It also
shows that every nontrivial connected graph in the class has ve-domination
number and isolation number one.

\begin{theorem}\label{thm:perfectly-characterization}
A graph $G$ is perfectly well-ve-dominated, or equivalently perfectly
well-isolated, if and only if $G$ is
$(P_4,K_1\vee 2K_2,\overline{P_3\cup P_3})$-free.
Furthermore, every nontrivial connected graph $G$ satisfying these
equivalent conditions has
\(
\gamma_{\mathrm{ve}}(G)=\iota(G)=1.
\)
\end{theorem}

\begin{proof}
The necessity follows immediately from
Corollary~\ref{cor:perfectly-forbidden}. For the converse, suppose that
$G$ is $(P_4,K_1\vee 2K_2,\overline{P_3\cup P_3})$-free. Since these
forbidden induced subgraph conditions are hereditary, it suffices to prove
that every such graph is well-ve-dominated. We proceed by induction on
$|V(G)|$. The assertion is immediate when $|V(G)|\leq 1$. Suppose now that
$|V(G)|\geq 2$ and that the assertion holds for all graphs of smaller order
satisfying the same forbidden induced subgraph conditions.
If $G$ is disconnected, then every nontrivial component of $G$ is a proper
induced subgraph satisfying the induction hypothesis, and hence is
well-ve-dominated. Proposition~\ref{prop:components} now implies that $G$ is
well-ve-dominated. Thus, we may assume that $G$ is connected.

Since $G$ is a connected $P_4$-free graph, it is a cograph. Therefore,
$G=G_1\vee G_2$, where $G_1$ and $G_2$ are nonempty induced subgraphs. We
first observe that at least one of $G_1$ and $G_2$ is $(K_2\cup K_1)$-free.
Otherwise, each would contain an induced $K_2\cup K_1$, and their join would
induce $(K_2\cup K_1)\vee(K_2\cup K_1)=\overline{P_3\cup P_3}$, contrary
to the assumption. Without loss of generality, assume that $G_1$ is
$(K_2\cup K_1)$-free.

We now show that every vertex of $G_1$ is a ve-dominating vertex of $G$.
Let $v\in V(G_1)$. Since $G=G_1\vee G_2$, the vertex $v$ is adjacent to
every vertex of $G_2$, and hence
$V(G)\setminus N_G[v]=V(G_1)\setminus N_{G_1}[v]$. If this set contained
an edge $xy$, then $\{v,x,y\}$ would induce a copy of $K_2\cup K_1$ in
$G_1$, a contradiction. Thus, $V(G)\setminus N_G[v]$ is independent.
Therefore, every edge of $G$ has at least one endpoint in $N_G[v]$, and so
$\{v\}$ is a ve-dominating set of $G$. Consequently, every minimal
ve-dominating set of $G$ intersecting $V(G_1)$ is a singleton.

Now let $D$ be a minimal ve-dominating set of $G$ such that
$D\cap V(G_1)=\varnothing$, and hence $D\subseteq V(G_2)$. Suppose first
that $G_2$ is edgeless. Then every vertex of $G_2$ is adjacent to every
vertex of $G_1$ and therefore is a ve-dominating vertex of $G$. Since $D$
is minimal, it follows that $|D|=1$.
We may therefore assume that $G_2$ contains an edge. We claim that $D$ is
a minimal ve-dominating set of $G_2$. Since $D\subseteq V(G_2)$ and $G_2$
is induced in $G$, every edge of $G_2$ that is ve-dominated by $D$ in $G$
is also ve-dominated by $D$ in $G_2$. Hence, $D$ ve-dominates $G_2$.
Suppose that a proper subset $D'\subsetneq D$ ve-dominates $G_2$. Since
$G_2$ contains an edge, $D'$ is nonempty. Every vertex of $D'$ is adjacent
to every vertex of $G_1$. Hence, $D'$ ve-dominates all edges between
$G_1$ and $G_2$, as well as all edges of $G_1$. Together with the fact
that $D'$ ve-dominates $G_2$, this implies that $D'$ ve-dominates $G$,
contradicting the minimality of $D$. Therefore, $D$ is a minimal
ve-dominating set of $G_2$.

Since $G$ is $(K_1\vee 2K_2)$-free and $G_1$ is nonempty, $G_2$ must be
$2K_2$-free; otherwise, an induced $2K_2$ in $G_2$ together with any
vertex of $G_1$ would induce a copy of $K_1\vee 2K_2$. Moreover, as an
induced subgraph of $G$, the graph $G_2$ is also
$(P_4,K_1\vee 2K_2,\overline{P_3\cup P_3})$-free.
If $G_2$ is connected, then
Lemma~\ref{lem:p4-2k2-free-vedom-vertex} implies that $G_2$ has a
ve-dominating vertex. If $G_2$ is disconnected, then, since it is
$2K_2$-free, it has at most one nontrivial component. This component is
connected, $P_4$-free, and $2K_2$-free, and hence has a ve-dominating
vertex by Lemma~\ref{lem:p4-2k2-free-vedom-vertex}; this vertex also
ve-dominates $G_2$. By the induction hypothesis, $G_2$ is
well-ve-dominated. Since $G_2$ has a ve-dominating vertex, it has a minimal
ve-dominating set of cardinality one. Therefore, every minimal
ve-dominating set of $G_2$ has cardinality one. In particular, $|D|=1$.

Consequently, every minimal ve-dominating set of $G$ has cardinality one.
Thus, $G$ is well-ve-dominated and, since $G$ is connected and nontrivial,
$\gamma_{\mathrm{ve}}(G)=1$. Since
$\gamma_{\mathrm{ve}}(G)=\iota(G)$, we also have $\iota(G)=1$.
Finally, every induced subgraph of $G$ satisfies the same forbidden induced
subgraph conditions and is therefore well-ve-dominated. Hence, $G$ is perfectly well-ve-dominated.
\end{proof}

\noindent
Theorem~\ref{thm:perfectly-characterization} provides a finite forbidden
induced subgraph characterization of perfectly well-ve-dominated graphs.
This is in sharp contrast to the class of well-ve-dominated graphs itself,
whose recognition problem is $\mathrm{co\text{-}NP}$-complete by
Theorem~\ref{thm:well-ve-dominated-complexity}. We next translate the two
non-$P_4$ obstructions into local conditions on a cotree. This gives both a
recursive structural description and a linear-time recognition algorithm.
For convenience, put $A=K_2\cup K_1$.

\begin{theorem}\label{thm:perfectly-cotree}
Let $G$ be a cograph and let $T_G$ be a reduced cotree of $G$. Then $G$ is
perfectly well-ve-dominated if and only if every join node $u$ of $T_G$
satisfies the following conditions:
\begin{enumerate}
    \item for every child $v$ of $u$, the graph $G_v$ is $2K_2$-free;
    \item at most one child $v$ of $u$ has the property that $G_v$ contains
    an induced copy of $A$.
\end{enumerate}
\end{theorem}

\begin{proof}
Suppose first that $G$ is perfectly well-ve-dominated, and let $u$ be a
join node of $T_G$. If a child $v$ of $u$ contained an induced $2K_2$,
then any vertex belonging to a different child of $u$, together with this
$2K_2$, would induce $K_1\vee 2K_2$, contradicting
Theorem~\ref{thm:perfectly-characterization}. Similarly, if two distinct
children of $u$ each contained an induced copy of $A$, then their union
would induce $A\vee A=(K_2\cup K_1)\vee(K_2\cup K_1)
=\overline{P_3\cup P_3}$, again a contradiction. Thus, both conditions
are necessary.

Conversely, suppose that every join node satisfies the two stated
conditions. Since $G$ is a cograph, it is $P_4$-free. Assume, for a
contradiction, that $G$ is not perfectly well-ve-dominated. By
Theorem~\ref{thm:perfectly-characterization}, $G$ contains an induced copy
of either $K_1\vee 2K_2$ or
$\overline{P_3\cup P_3}=A\vee A$.
First suppose that $G$ contains an induced $K_1\vee 2K_2$, and let $u$ be
the lowest cotree node whose subtree contains all five vertices of this
copy. Since $K_1\vee 2K_2$ is connected, $u$ is a join node. The four
vertices inducing the $2K_2$ must lie in a single child of $u$, since
$2K_2$ has no nontrivial join decomposition: its complement $C_4$ is
connected. By the minimal choice of $u$, the universal vertex cannot lie
in that same child. Hence, some child of $u$ contains an induced $2K_2$,
contradicting condition~(1).
Now suppose that $G$ contains an induced $A\vee A$, and let $u$ be the
lowest cotree node whose subtree contains all six vertices of this copy.
Again, $u$ is a join node. Since $\overline{A}=P_3$ is connected, $A$ has
no nontrivial join decomposition, so the vertices of each induced copy of
$A$ lie entirely within a single child of $u$. By the minimal choice of
$u$, the two copies cannot lie in the same child and therefore lie in
distinct children. Thus, two children of $u$ contain induced copies of
$A$, contradicting condition~(2).
Hence, neither obstruction occurs. By
Theorem~\ref{thm:perfectly-characterization}, $G$ is perfectly
well-ve-dominated.
\end{proof}

For an internal cotree node $u$, let $C(u)$ denote the set of its children.
Let $e(u)$, $a(u)$, and $b(u)$ indicate, respectively, whether $G_u$
contains an edge, an induced copy of $A$, and an induced copy of $2K_2$.
These three Boolean values can be computed in a single postorder traversal
of the cotree.

\begin{algorithm}[ht]
\caption{\textsc{Perfectly Well-VE-Dominated Recognition}}
\label{alg:perfectly-recognition}
\begin{algorithmic}[1]
\Require A graph $G$
\Ensure \textsc{Yes} if $G$ is perfectly well-ve-dominated, and
\textsc{No} otherwise

\State Attempt to construct a reduced cotree $T_G$ of $G$
\If{no such cotree exists}
    \State \Return \textsc{No}
\EndIf

\For{each node $u$ of $T_G$ in postorder}
    \If{$u$ is a leaf}
        \State $e(u)\gets\textsc{False}$,
        $a(u)\gets\textsc{False}$,
        $b(u)\gets\textsc{False}$
    \ElsIf{$u$ is a union node with children $C(u)$}
        \State $e(u)\gets \bigvee_{v\in C(u)}e(v)$
        \State $a(u)\gets
        \left(\bigvee_{v\in C(u)}a(v)\right)\vee e(u)$
        \State $b(u)\gets
        \left(\bigvee_{v\in C(u)}b(v)\right)\vee
        \left(|\{v\in C(u):e(v)\}|\geq2\right)$
    \Else
        \If{$\bigvee_{v\in C(u)}b(v)$ or
        $|\{v\in C(u):a(v)\}|\geq2$}
            \State \Return \textsc{No}
        \EndIf
        \State $e(u)\gets\textsc{True}$
        \State $a(u)\gets\bigvee_{v\in C(u)}a(v)$
        \State $b(u)\gets\bigvee_{v\in C(u)}b(v)$
    \EndIf
\EndFor

\State \Return \textsc{Yes}
\end{algorithmic}
\end{algorithm}

\begin{proposition}\label{prop:perfectly-recognition}
Algorithm~\ref{alg:perfectly-recognition} correctly recognizes perfectly
well-ve-dominated graphs in time $O(n+m)$, where $n=|V(G)|$ and
$m=|E(G)|$.
\end{proposition}

\begin{proof}
We first verify by induction in postorder that $e(u)$, $a(u)$, and $b(u)$
have their stated meanings. This is immediate at a leaf.
Suppose that $u$ is a union node. An induced copy of $A=K_2\cup K_1$
either lies entirely within one child or consists of an edge in one child
and a vertex in another. Since $T_G$ is reduced, every internal node has at
least two children; hence, whenever $G_u$ contains an edge, such a vertex
is available in another child. Thus,
$a(u)=\big(\bigvee_{v\in C(u)}a(v)\big)\vee e(u)$. Similarly, an induced
$2K_2$ either lies within one child or consists of an edge from each of two
distinct children, which gives the stated update for $b(u)$. The update for
$e(u)$ is immediate.
Now suppose that $u$ is a join node. Since $\overline{A}=P_3$ and
$\overline{2K_2}=C_4$ are connected, neither $A$ nor $2K_2$ has a
nontrivial join decomposition. Hence, any induced copy of either graph in
$G_u$ lies entirely within one child, so $a(u)$ and $b(u)$ are precisely
the disjunctions of the corresponding child values. Moreover, since a
reduced cotree has at least two children at every internal node, a join
node necessarily contains an edge, and hence $e(u)=\textsc{True}$.

If the input graph is not a cograph, then it contains an induced $P_4$ and
is not perfectly well-ve-dominated by
Theorem~\ref{thm:perfectly-characterization}. Otherwise, at each join node
the algorithm rejects precisely when some child contains an induced
$2K_2$ or at least two children contain induced copies of $A$. By
Theorem~\ref{thm:perfectly-cotree}, these are exactly the obstructions to
$G$ being perfectly well-ve-dominated. Thus, the algorithm returns
\textsc{Yes} if and only if $G$ is perfectly well-ve-dominated.

A reduced cotree can be constructed in $O(n+m)$ time. It has $O(n)$ nodes
and $O(n)$ child links, and all Boolean values and join-node tests can be
computed by scanning the children of each node once. Hence, the postorder
phase takes $O(n)$ time, and the overall running time is $O(n+m)$.
\end{proof}

\begin{corollary}\label{cor:perfectly-linear}
Perfectly well-ve-dominated graphs can be recognized in linear time.
\end{corollary}

A direct application of
Theorem~\ref{thm:perfectly-characterization} provides a structural
characterization of triangle-free perfectly well-ve-dominated graphs.

\begin{corollary}\label{cor:triangle-free-characterization}
A connected nontrivial triangle-free graph is perfectly
well-ve-dominated if and only if it is complete bipartite.
\end{corollary}

\begin{proof}
Let $G$ be a connected nontrivial triangle-free perfectly
well-ve-dominated graph. By
Theorem~\ref{thm:perfectly-characterization}, $G$ is $P_4$-free, and hence
$G$ is a connected cograph. Therefore, $G=G_1\vee G_2$ for some nonempty
induced subgraphs $G_1$ and $G_2$.
Since $G$ is triangle-free, both $G_1$ and $G_2$ must be edgeless. Indeed, if
one of them contained an edge, then that edge together with any vertex of the
other subgraph would induce a triangle in $G$. Hence,
$G=\overline{K_{|V(G_1)|}}\vee\overline{K_{|V(G_2)|}}$.
Therefore, $G$ is complete bipartite.
Conversely, every complete bipartite graph contains no induced
$P_4$, $K_1\vee 2K_2$, or $\overline{P_3\cup P_3}$. Indeed, an induced
subgraph of a complete bipartite graph is complete bipartite and hence cannot
be a $P_4$, while each of the latter two forbidden graphs contains a
triangle. Therefore, by
Theorem~\ref{thm:perfectly-characterization}, every complete bipartite graph
is perfectly well-ve-dominated.
\end{proof}

\noindent
Corollary~\ref{cor:triangle-free-characterization} also yields a particularly
simple direct recognition procedure for the connected triangle-free case.

\begin{corollary}
Connected nontrivial triangle-free perfectly well-ve-dominated graphs can
be recognized in linear time.
\end{corollary}

\begin{proof}
By Corollary~\ref{cor:triangle-free-characterization}, it suffices to test
whether the input graph is complete bipartite.
A breadth-first search determines whether the graph is bipartite and
produces a bipartition $(A,B)$ in time $O(n+m)$, where $n=|V(G)|$ and
$m=|E(G)|$. The graph is complete
bipartite if and only if
\(
|E(G)|=|A||B|.
\)
Hence the recognition algorithm runs in $O(n+m)$ time.
\end{proof}

\section{Well-ve-Dominated Trees}\label{sec:WEDtrees}

We next turn to trees, where well-ve-domination admits a complete structural characterization. Since well-ve-domination and well-isolation
are equivalent, we state the results in terms of well-ve-dominated trees;
the corresponding well-isolated formulations follow immediately. We also
derive a canonical representation of every reduced well-ve-dominated tree
and obtain a linear-time recognition algorithm.

We first recall the tree characterization of well-ve-covered graphs from
\cite{BoutrigChellaliMeddah2025}. Let $\mathcal{T}_{\mathrm{wc}}$ be the
family obtained from $r\geq1$ pairwise vertex-disjoint stars, each of order
at least three, by adding $r-1$ edges between leaves of distinct stars so
that the resulting graph is a tree and, for each original star, all added
edges incident to that star meet the same selected leaf.
Figure~\ref{fig:family-F} illustrates a member of
$\mathcal{T}_{\mathrm{wc}}$ obtained from six stars.

\begin{theorem}[{\cite{BoutrigChellaliMeddah2025}}]
\label{thm:well-vecovered-trees}
A nontrivial tree $T$ is well-ve-covered if and only if $T=P_2$ or
$T\in\mathcal{T}_{\mathrm{wc}}$.
\end{theorem}

\noindent
Since every well-ve-dominated graph is necessarily well-ve-covered,
Theorem~\ref{thm:well-vecovered-trees} implies that every nontrivial
well-ve-dominated tree belongs to
$\mathcal{T}_{\mathrm{wc}}\cup\{P_2\}$.
The tree $P_2$ is trivially well-ve-dominated. Therefore, to prove the
converse and establish the coincidence of well-ve-covered and
well-ve-dominated trees, it remains to show that every member of
$\mathcal{T}_{\mathrm{wc}}$ is well-ve-dominated.

\begin{figure}[t]
\centering
\begin{tikzpicture}[
    scale=0.75,
    transform shape,
    every node/.style={font=\small},
    center/.style={
        circle,
        draw,
        fill=blue!10,
        minimum size=4.5mm,
        inner sep=0pt
    },
    leaf/.style={
        circle,
        draw,
        fill=green!12,
        minimum size=4.5mm,
        inner sep=0pt
    },
    connector/.style={
        circle,
        draw,
        fill=orange!18,
        minimum size=4.5mm,
        inner sep=0pt
    },
    edge/.style={
        line width=0.9pt
    },
    addededge/.style={
        line width=1.1pt,
        dashed
    }
]

% Star 1
\node[center]    (c1) at (0,0) {$c_1$};
\node[leaf]      (a1) at (-1.0,1.0) {};
\node[leaf]      (a2) at (-1.2,0) {};
\node[leaf]      (a3) at (-1.0,-1.0) {};
\node[connector] (a4) at (1.0,0) {};

\draw[edge] (c1)--(a1);
\draw[edge] (c1)--(a2);
\draw[edge] (c1)--(a3);
\draw[edge] (c1)--(a4);

% Star 2
\node[center]    (c2) at (3.7,0) {$c_2$};
\node[connector] (b1) at (2.7,0) {};
\node[leaf]      (b2) at (3.7,1.2) {};
\node[leaf]      (b3) at (4.7,1.0) {};
\node[leaf]      (b4) at (3.7,-1.2) {};
\node[leaf]      (b5) at (4.7,-1.0) {};

\draw[edge] (c2)--(b1);
\draw[edge] (c2)--(b2);
\draw[edge] (c2)--(b3);
\draw[edge] (c2)--(b4);
\draw[edge] (c2)--(b5);

% Star 3
\node[center]    (c3) at (0,-4) {$c_3$};
\node[leaf]      (d1) at (-1,-4) {};
\node[connector] (d2) at (0,-2.8) {};
\node[leaf]      (d3) at (1,-3) {};
\node[leaf]      (d4) at (0,-5.2) {};
\node[leaf]      (d5) at (1,-5) {};

\draw[edge] (c3)--(d1);
\draw[edge] (c3)--(d2);
\draw[edge] (c3)--(d3);
\draw[edge] (c3)--(d4);
\draw[edge] (c3)--(d5);

% Star 4
\node[center]    (c4) at (3.7,-4.0) {$c_4$};
\node[connector] (e1) at (2.7,-4.0) {};
\node[leaf]      (e2) at (3.7,-2.8) {};
\node[leaf]      (e3) at (4.7,-3.0) {};
\node[leaf]      (e4) at (3.7,-5.2) {};
\node[leaf]      (e5) at (4.7,-5.0) {};

\draw[edge] (c4)--(e1);
\draw[edge] (c4)--(e2);
\draw[edge] (c4)--(e3);
\draw[edge] (c4)--(e4);
\draw[edge] (c4)--(e5);

% Star 5
\node[center]    (c5) at (-3.5,-4.0) {$c_5$};
\node[connector] (f1) at (-2.5,-4.0) {};
\node[leaf]      (f2) at (-3.5,-2.8) {};
\node[leaf]      (f3) at (-4.5,-3.0) {};
\node[leaf]      (f4) at (-3.5,-5.2) {};
\node[leaf]      (f5) at (-2.5,-5.0) {};

\draw[edge] (c5)--(f1);
\draw[edge] (c5)--(f2);
\draw[edge] (c5)--(f3);
\draw[edge] (c5)--(f4);
\draw[edge] (c5)--(f5);

% Star 6
\node[center]    (c6) at (-4.5,0) {$c_6$};
\node[leaf]      (g1) at (-3.5,0) {};
\node[leaf]      (g2) at (-4.5,1.2) {};
\node[leaf]      (g3) at (-5.5,-1) {};
\node[leaf]      (g4) at (-4.5,-1.2) {};
\node[connector] (g5) at (-3.5,-1) {};

\draw[edge] (c6)--(g1);
\draw[edge] (c6)--(g2);
\draw[edge] (c6)--(g3);
\draw[edge] (c6)--(g4);
\draw[edge] (c6)--(g5);

% Added edges between selected leaves
\draw[addededge] (a4)--(b1);
\draw[addededge] (b1)--(d2);
\draw[addededge] (b1)--(e1);
\draw[addededge] (d2)--(f1);
\draw[addededge] (g5)--(f1);

\end{tikzpicture}

\caption{A tree in the family $\mathcal{T}_{\mathrm{wc}}$ obtained from six pairwise
vertex-disjoint stars.}
\label{fig:family-F}
\end{figure}

\begin{proposition}\label{prop:Twc-well-ve-dominated}
Every tree in the family $\mathcal{T}_{\mathrm{wc}}$ is
well-ve-dominated.
\end{proposition}

\begin{proof}
Let $T\in\mathcal{T}_{\mathrm{wc}}$. Then $T$ is obtained from $r$ pairwise
vertex-disjoint stars, each of order at least three, by adding $r-1$ edges
between selected leaves so that the resulting graph is a tree. For each
original star, all added edges incident to that star are incident to the
same selected leaf.
If $r=1$, then $T$ is a star. Its reduction is $P_2$, which is
well-ve-dominated. Hence, by Lemma~\ref{lem:reduced}, the tree $T$ is
well-ve-dominated. We may therefore assume that $r\geq 2$.

By Lemma~\ref{lem:reduced}, it is enough to consider the reduced tree
$T'=R(T)$. For each original star, let $c_i$ be its center, let $w_i$ be the
selected leaf used for the added edges, and let $\ell_i$ be the unique
remaining leaf adjacent to $c_i$ in $T'$. Then
\[
V(T')=\{\ell_i,c_i,w_i:1\leq i\leq r\} \quad
\text{ and } \quad
E(T')
=
\{\ell_i c_i,c_iw_i:1\leq i\leq r\}
\cup E(H),
\]
where $H=T'[\{w_1,\ldots,w_r\}]$ is a tree. Thus, $T'$ consists of the
backbone tree $H$ together with a pendant path
\(
\ell_i c_i w_i
\)
attached to each backbone vertex $w_i$. In particular, each edge
$\ell_i c_i$ is a good pendant edge of $T'$.

Let $D$ be a minimal ve-dominating set of $T'$. Since the pendant edge
$\ell_i c_i$ can be ve-dominated only by vertices in
\(
\{\ell_i,c_i,w_i\},
\)
we have
\(
D\cap\{\ell_i,c_i,w_i\}\neq\emptyset
\)
for every $i$. Hence, $|D|\geq r$.
We next show that $D$ contains at most one vertex from each set
$\{\ell_i,c_i,w_i\}$. If $\ell_i\in D$ together with either $c_i$ or $w_i$,
then $\ell_i$ has no private edge with respect to $D$. Similarly, if both
$c_i$ and $w_i$ belong to $D$, then $c_i$ has no private edge, since every
edge ve-dominated by $c_i$ is also ve-dominated by $w_i$. In either case,
Proposition~\ref{prop:minimal-ved} contradicts the minimality of $D$.
Therefore,
\(
|D\cap\{\ell_i,c_i,w_i\}|=1
\)
for every $i$, and consequently $|D|=r$.

Therefore, every minimal ve-dominating set of $T'$ has cardinality $r$, and hence
$T'$ is well-ve-dominated. By Lemma~\ref{lem:reduced}, the original tree
$T$ is also well-ve-dominated.
\end{proof}

\noindent
Combining Proposition~\ref{prop:Twc-well-ve-dominated} with
Theorem~\ref{thm:well-vecovered-trees} gives the main characterization of
this section.

\begin{theorem}\label{thm:well-ve-dominated-trees}
A tree $T$ is well-ve-dominated if and only if
$T\in\mathcal{T}_{\mathrm{wc}}\cup\{P_1,P_2\}$.
\end{theorem}

\begin{proof}
Suppose first that $T$ is well-ve-dominated. If $T=P_1$, there is nothing
to prove. Otherwise, every well-ve-dominated graph is well-ve-covered, and
hence Theorem~\ref{thm:well-vecovered-trees} implies that $T=P_2$ or
$T\in\mathcal{T}_{\mathrm{wc}}$.
Conversely, the trees $P_1$ and $P_2$ are well-ve-dominated, while every
tree in $\mathcal{T}_{\mathrm{wc}}$ is well-ve-dominated by
Proposition~\ref{prop:Twc-well-ve-dominated}. Therefore, the stated
characterization follows.
\end{proof}

\begin{corollary}\label{cor:well-ve-paths}
A path $P_n$ is well-ve-dominated if and only if
$n\in\{1,2,3,6\}$.
\end{corollary}

\begin{proof}
By Theorem~\ref{thm:well-ve-dominated-trees}, it remains to determine which
paths belong to $\mathcal{T}_{\mathrm{wc}}$. Since every vertex of a path
has degree at most two, each original star in the construction of
$\mathcal{T}_{\mathrm{wc}}$ must have order three. Moreover, the construction
cannot involve three or more stars: the tree induced by the selected leaves
would then have a vertex of degree at least two, and the corresponding
selected leaf would have degree at least three in the resulting tree.
Therefore, the construction uses either one star, yielding $P_3$, or two
stars, yielding $P_6$. Together with the exceptional paths $P_1$ and $P_2$,
the result follows.
\end{proof}

\begin{corollary}\label{cor:well-covered-well-dominated-trees}
Among nontrivial trees, well-ve-coveredness and well-ve-domination
coincide.
\end{corollary}

The proof of Proposition~\ref{prop:Twc-well-ve-dominated} also reveals a
particularly simple canonical form for reduced well-ve-dominated trees.
After applying the reduction described in Lemma~\ref{lem:reduced}, the
selected leaves used to connect the original stars form a backbone tree.
Each backbone vertex $w_i$ is incident to a unique pendant path
$\ell_i c_i w_i$,
where $\ell_i c_i$ is a good pendant edge.
Figure~\ref{fig:reduced-F-backbone} illustrates this canonical structure.

\begin{figure}[b]
\centering
\begin{tikzpicture}[
    scale=0.75,
    transform shape,
    every node/.style={font=\small},
    wvertex/.style={
        circle,
        draw,
        fill=blue!10,
        minimum size=5mm,
        inner sep=0pt
    },
    center/.style={
        circle,
        draw,
        fill=green!12,
        minimum size=5mm,
        inner sep=0pt
    },
    leaf/.style={
        circle,
        draw,
        fill=green!12,
        minimum size=6mm,
        inner sep=0pt
    },
    edge/.style={
        line width=0.9pt
    }
]

% Backbone vertices
\node[wvertex] (w1) at (0,0) {$w_1$};
\node[wvertex] (w2) at (1.8,0) {$w_2$};
\node[wvertex] (w3) at (3.6,0) {$w_3$};
\node[wvertex] (w4) at (5.4,0) {$w_4$};
\node[wvertex] (w5) at (7.2,0) {$w_5$};
\node[wvertex] (w6) at (3.6,-1.5) {$w_6$};
\node[wvertex] (w7) at (5.4,-1.5) {$w_7$};

% Backbone edges
\draw[edge] (w1)--(w2)--(w3)--(w4)--(w5);
\draw[edge] (w3)--(w6)--(w7);

% Pendant paths
\node[center] (c1) at (0,0.9) {$c_1$};
\node[leaf]   (l1) at (0,1.8) {$\ell_1$};
\draw[edge] (w1)--(c1)--(l1);

\node[center] (c2) at (1.8,0.9) {$c_2$};
\node[leaf]   (l2) at (1.8,1.8) {$\ell_2$};
\draw[edge] (w2)--(c2)--(l2);

\node[center] (c3) at (3.6,0.9) {$c_3$};
\node[leaf]   (l3) at (3.6,1.8) {$\ell_3$};
\draw[edge] (w3)--(c3)--(l3);

\node[center] (c4) at (5.4,0.9) {$c_4$};
\node[leaf]   (l4) at (5.4,1.8) {$\ell_4$};
\draw[edge] (w4)--(c4)--(l4);

\node[center] (c5) at (8.1,0) {$c_5$};
\node[leaf]   (l5) at (9.0,0) {$\ell_5$};
\draw[edge] (w5)--(c5)--(l5);

\node[center] (c6) at (3.6,-2.4) {$c_6$};
\node[leaf]   (l6) at (3.6,-3.3) {$\ell_6$};
\draw[edge] (w6)--(c6)--(l6);

\node[center] (c7) at (6.3,-1.5) {$c_7$};
\node[leaf]   (l7) at (7.2,-1.5) {$\ell_7$};
\draw[edge] (w7)--(c7)--(l7);

\end{tikzpicture}

\caption{The canonical representation of a reduced tree in
$\mathcal{T}_{\mathrm{wc}}$.}
\label{fig:reduced-F-backbone}
\end{figure}

\begin{theorem}\label{thm:reduced-well-ve-dominated-trees}
For a reduced tree $T$, $T$ is well-ve-dominated if and only if one
of the following holds:
\begin{enumerate}
    \item $T\in\{P_1,P_2\}$;
    \item $T$ is obtained from a tree $H$ of order at least two by attaching
    a pendant path of length two to each vertex of $H$.
\end{enumerate}
\end{theorem}

\begin{proof}
Suppose first that $T$ is a reduced well-ve-dominated tree. If
$|V(T)|\leq 2$, then $T\in\{P_1,P_2\}$. Assume that $|V(T)|\geq 3$.
By Theorem~\ref{thm:well-ve-dominated-trees},
$T\in\mathcal{T}_{\mathrm{wc}}$. Since $T$ is reduced, each original
star contributes exactly three vertices: its center, the selected leaf used
to connect it to the other stars, and one additional leaf. The selected
leaves induce a tree $H$, and each of its vertices is incident to a unique
pendant path of length two.
Conversely, suppose that $T$ is obtained from a tree $H$ by attaching a
pendant path of length two to each vertex of $H$. Then
$T\in\mathcal{T}_{\mathrm{wc}}$, and hence $T$ is well-ve-dominated by
Proposition~\ref{prop:Twc-well-ve-dominated}. The trees $P_1$ and $P_2$
are also well-ve-dominated.
\end{proof}

Let $T$ be a reduced well-ve-dominated tree of order at least three.
Its canonical representation yields a natural partition of its vertex set.
Let $L_T$ denote the set of leaves incident to good pendant edges, let
$S_T$ denote their support vertices, and define
$W_T=N_T(S_T)\setminus L_T$. Then $W_T$ is precisely the vertex set of
the backbone tree. Moreover, each vertex of $S_T$ has exactly one neighbor
in $L_T$ and exactly one neighbor in $W_T$, while each vertex of $W_T$
has exactly one neighbor in $S_T$.

\begin{corollary}\label{cor:structure}
Let $T$ be a reduced well-ve-dominated tree.
\begin{enumerate}
    \item If $|V(T)|\leq 2$, then
    $\gamma_{\mathrm{ve}}(T)=\iota(T)=0$ when $T=P_1$, while
    $\gamma_{\mathrm{ve}}(T)=\iota(T)=1$ when $T=P_2$.

    \item If $|V(T)|\geq 3$, then
    $|L_T|=|S_T|=|W_T|=\gamma_{\mathrm{ve}}(T)=\iota(T)$, and thus
    $|V(T)|=3\gamma_{\mathrm{ve}}(T)=3\iota(T)$.
\end{enumerate}
\end{corollary}

The canonical partition also gives an explicit realization of a result of
Boyer and Goddard~\cite{BoyerGoddard2024}, who proved that every connected
graph of order at least six has three pairwise disjoint isolating sets. In
the present setting, these sets form a partition of the vertex set and are
all minimum isolating sets.

\begin{corollary}\label{cor:three-isolating-sets}
Let $T$ be a reduced well-ve-dominated tree of order at least three.
Let $(X,Y)$ be the bipartition of the backbone tree $T[W_T]$, and define
$L_X=\{\ell_i:w_i\in X\}$ and $L_Y=\{\ell_i:w_i\in Y\}$. Then
$S_T$, $X\cup L_Y$, and $Y\cup L_X$ are three pairwise disjoint minimum
isolating sets whose union is $V(T)$.
\end{corollary}

\begin{proof}
The three sets are pairwise disjoint and their union is $V(T)$. The set
$S_T$ is a dominating set and hence an isolating set. For $X\cup L_Y$,
every support vertex is adjacent either to a vertex of $X$ or to its leaf
in $L_Y$, while every vertex of $Y$ has a neighbor in $X$, since the
backbone is a nontrivial tree. Hence,
$V(T)\setminus N_T[X\cup L_Y]\subseteq L_X$, which is independent.
Thus, $X\cup L_Y$ is isolating. By symmetry, $Y\cup L_X$ is also
isolating. Finally,
$|S_T|=|X\cup L_Y|=|Y\cup L_X|=|W_T|=\iota(T)$ by
Corollary~\ref{cor:structure}. Therefore, all three sets are minimum
isolating sets.
\end{proof}

We now make the connection with extremal isolation explicit. Let
$\mathcal{T}$ be the family of trees obtained from an arbitrary tree $T_0$
by attaching a pendant path of length two to each vertex of $T_0$;
equivalently, for every $v\in V(T_0)$, add two vertices $s_v$ and $\ell_v$
and the edges $vs_v$ and $s_v\ell_v$. The same family can be obtained by
starting with a path $\ell s v\cong P_3$ and repeatedly adding a new path
$\ell's'v'\cong P_3$ together with an edge $vv'$, where $v$ is an existing
backbone vertex and $v'$ becomes a new backbone vertex. Thus, the
constructive formulations used in the two cited characterizations define
the same family.
We use the following characterization of the trees attaining equality in
the sharp isolation bound.

\begin{theorem}[{\cite{krishnakumari2014,LemanskaEtAl2021}}]
\label{thm:extremal-isolation-trees}
Let $T$ be a tree of order $n\geq 3$. Then $\iota(T)=n/3$ if and only if
$T\in\mathcal{T}$.
\end{theorem}

\begin{corollary}\label{cor:reduced-extremal-isolation}
For a reduced tree $T$ of order at least three, the following
statements are equivalent:
\begin{enumerate}
    \item $T$ is well-ve-dominated;
    \item $T\in\mathcal{T}$;
    \item $\iota(T)=|V(T)|/3$.
\end{enumerate}
\end{corollary}

\begin{proof}
By Theorem~\ref{thm:reduced-well-ve-dominated-trees}, a reduced tree of
order at least three is well-ve-dominated if and only if it is obtained by
attaching a pendant path of length two to every vertex of a backbone tree.
This is precisely the construction defining $\mathcal{T}$, so $(1)$ and
$(2)$ are equivalent. The equivalence of $(2)$ and $(3)$ follows from
Theorem~\ref{thm:extremal-isolation-trees}.
\end{proof}

\noindent
Corollary~\ref{cor:reduced-extremal-isolation} shows that reduced
well-ve-dominated trees of order at least three are precisely the extremal
trees attaining equality in the sharp isolation bound
$\iota(T)\leq |V(T)|/3$~\cite{CaroHansberg2017,Zylinski2019}.

The canonical form also yields a direct recognition procedure. After
reducing the input tree, it remains only to determine whether its vertices
admit the partition $L_T\cup S_T\cup W_T$ described above, with $W_T$
inducing the backbone tree.

\begin{corollary}\label{cor:linear-recognition-trees}
Given a tree $T$, it can be decided in linear time whether $T$ is
well-ve-dominated.
\end{corollary}

\begin{proof}
In a tree of order at least three, two vertices are false twins if and only
if they are leaves adjacent to the same support vertex. Indeed, false twins
are nonadjacent, and if they had two common neighbors, those four vertices
would lie on a cycle; hence their common open neighborhood consists of a
single vertex. The converse is immediate. Thus, the reduced tree
$T'=R(T)$ is obtained by deleting all but one leaf adjacent to each support
vertex. This can be done in linear time by inspecting the degrees and
adjacency lists.

If $T'\cong P_1$ or $T'\cong P_2$, then $T$ is well-ve-dominated. Assume
therefore that $|V(T')|\geq 3$. Let $L$ be the set of leaves of $T'$, let
$S=N_{T'}(L)$, and define $W=N_{T'}(S)\setminus L$. Since $T'$ is
reduced, each vertex of $S$ is adjacent to exactly one leaf. We check
whether $|W|\geq 2$, whether $V(T')=L\cup S\cup W$ is a partition, whether
every vertex of $S$ has exactly one neighbor in $L$ and exactly one neighbor
in $W$, whether every vertex of $W$ has exactly one neighbor in $S$, and
whether $T'[W]$ is a tree. All these conditions can be verified in linear
time by traversing $T'$.

If one of these conditions fails, then $T'$ does not have the canonical
form of Theorem~\ref{thm:reduced-well-ve-dominated-trees} and hence is not
well-ve-dominated. If all conditions hold, then $T'$ consists of the
backbone tree $T'[W]$ together with one pendant path $\ell c w$ attached
to each backbone vertex $w\in W$, and hence is well-ve-dominated.
Finally, Lemma~\ref{lem:reduced} implies that $T$ is well-ve-dominated if
and only if $T'$ is well-ve-dominated. Therefore, the recognition procedure
runs in linear time.
\end{proof}

\section{Conclusion}

In this paper, we introduce well-ve-dominated graphs and developed their
structural, algorithmic, and isolation-theoretic foundations. Equivalently,
these are the graphs whose minimal isolating sets all have the same
cardinality, providing a natural uniformity counterpart to the classical
theory of graph isolation.
We show that the general recognition problem is $\mathrm{co\text{-}NP}$-complete; in
contrast, the largest hereditary subclass of well-ve-dominated graphs admits
an exact characterization by three forbidden induced subgraphs and a
linear-time recognition algorithm. Every connected nontrivial graph in this
subclass has vertex-edge domination and isolation number one, while its
triangle-free members are precisely the complete bipartite graphs.

For trees, we show that well-ve-coveredness and well-ve-domination coincide on every
nontrivial tree. After false-twin reduction, every well-ve-dominated tree of
order at least three has a canonical form obtained by attaching a pendant
path of length two to each vertex of a backbone tree. This representation
yields a linear-time recognition algorithm and identifies the reduced
well-ve-dominated trees precisely with the extremal family satisfying
$\iota(T)=|V(T)|/3$. It also yields a partition into three pairwise disjoint
minimum isolating sets. Thus, in the reduced tree setting, uniformity among
minimal isolating sets and extremality for the isolation number describe
exactly the same family.

Therefore, these results exhibit a sharp contrast: well-ve-domination
is computationally difficult to recognize in general, yet imposes strong
structure under hereditary or acyclic restrictions. This suggests several
directions for further study, including characterizations and the complexity of recognition on other natural
graph classes, and the
identification of further classes on which well-ve-coveredness and
well-ve-domination coincide.

\section*{Acknowledgements}
The author sincerely thanks Professor Wayne Goddard for his constructive
comments and Professor Didem Gözüpek for her valuable suggestions and
insightful discussions, all of which greatly improved this work.

\section*{Funding}
This work is supported by the Scientific and Technological Research Council of Turkey (TUBITAK) under grant no. 126F284.

\section*{Data availability}
No data were generated or analyzed during the current study.

\section*{Declaration of competing interest}
The author declares that she has no conflict of interest.

\bibliographystyle{elsarticle-num}
\bibliography{References_submission}

\end{document}